\documentclass{article}

\usepackage{amsmath}
\usepackage{amssymb}
\usepackage{theorem}
\usepackage{xypic}
\xyoption{all}
\usepackage{enumitem}
\usepackage{rotating}
\usepackage{mathrsfs}
\usepackage{stmaryrd}
\usepackage{tikz}

\newenvironment{items}
{\begin{enumerate}[topsep=3pt, itemsep=3pt, parsep=0pt, label=(\roman*)]}
{\end{enumerate}}
\renewcommand{\tilde}{\widetilde}

\newcommand{\XX}{{\mathfrak X}}

\newcommand{\UU}{{\mathfrak U}}

\newcommand{\Z}{{\mathscr Z}}
\newcommand{\E}{{\mathscr E}}
\newcommand{\F}{{\mathscr F}}

\newcommand{\LL}{{\mathscr L}}
\newcommand{\R}{{\mathscr R}}

\newcommand{\MM}{{\mathfrak M}}

\newcommand{\m}{{\mathfrak m}}

\newcommand{\zz}{{\mathbb Z}}

\newcommand{\qq}{{\mathbb Q}}

\newcommand{\cc}{{\mathbb C}}

\newcommand{\OO}{{\mathscr O}}

\newcommand{\gG}{{\mathscr G}}
\newcommand{\fF}{{\mathscr F}}
\newcommand{\lL}{{\mathscr L}}

\newcommand{\gr}{{\mathop{\rm gr}\nolimits}}
\renewcommand{\sp}{{\mathop{\rm sp}\nolimits}}

\newcommand{\Mor}{\mathop{\rm Mor}\nolimits}

\newcommand{\sheafhom}{\mathop{{\mathscr H\! om}}\nolimits}
\newcommand{\sheafend}{\mathop{{\mathscr E\! nd}}\nolimits}

\newcommand{\End}{\mathop{\rm End}\nolimits}
\newcommand{\GL}{\mathop{\rm GL}\nolimits}

\newcommand{\projectlim}{\mathop{{\lim\limits_{\textstyle\longleftarrow}}}\limits}

\newcommand{\hooklongrightarrow}{\lhook\joinrel\longrightarrow}

\newtheorem{thm}{Theorem}[section]
\newtheorem{cor}[thm]{Corollary}
\newtheorem{lem}[thm]{Lemma}
\newtheorem{prop}[thm]{Proposition}

\theorembodyfont{\upshape}
\newtheorem{defn}[thm]{Definition}

\newtheorem{rmk}[thm]{Remark}

\newtheorem{ex}[thm]{Example}

\newenvironment{pf}{\begin{trivlist}\item[]{\sc Proof.}}%
            {\nolinebreak $\Box$ \end{trivlist}}

            {\nolinebreak $\Box$ \end{trivlist}}

\newcommand{\com}{\bullet}

\newcommand{\Proj}{\mathop{\mathbb{P}\rm{roj}}}

\newcommand{\rk}{\mathop{\rm rk}\nolimits}

\newcommand{\spec}{\mathop{\rm Spec}\nolimits}

\newcommand{\Sym}{\mathop{\rm Sym}\nolimits}
\newcommand{\Ext}{\mathop{\rm Ext}\nolimits}
\newcommand{\Hom}{\mathop{\rm Hom}\nolimits}

\newcommand{\sfrac}[2]{{\textstyle\frac{#1}{#2}}}

\newcommand{\noprint}[1]{}

\def\Label#1{\label{#1}{\tt [#1]}\phantom{Mh}}

\def\Label{\label}

\newcommand{\dg}{differential graded }
\newcommand{\RAct}{\mathop{{\frak R}{\rm Act}}\nolimits}
\newcommand{\Act}{\mathop{{\rm Act}}\nolimits}
\newcommand{\RMod}{\mathop{{\frak R}{\frak M}{\rm od}}\nolimits}
\newcommand{\Mod}{\mathop{{\frak M}{\rm od}}\nolimits}
\newcommand{\Modtilde}{\mathop{\tilde{\frak M}{\rm od}}\nolimits}
\newcommand{\RModtilde}{\mathop{{\frak R}{\tilde{\frak M}}{\rm od}}\nolimits}
\newcommand{\RQuot}{\mathop{{\frak R}{\rm Quot}}\nolimits}

\newcommand{\git}{\!\sslash\!}

\author{K. Behrend, I. Ciocan-Fontanine, J. Hwang, M. Rose}
\title{The Derived Moduli Space of Stable Sheaves}
\date{March 31, 2010}

\begin{document}
\sloppy

\maketitle

\begin{abstract} {We construct the derived scheme of stable sheaves on
    a smooth projective variety via derived moduli of finite graded
    modules over a graded ring. We do this by dividing the derived
    scheme of actions of Ciocan-Fontanine and Kapranov by a suitable
    algebraic gauge group. We show that the natural notion of
    GIT-stability for graded modules reproduces stability for
    sheaves.}
\end{abstract}

\tableofcontents

\section*{Introduction}
\addcontentsline{toc}{section}{Introduction}

For some years it has been a tenant of geometry, that
deformation theory problems are governed by differential graded Lie
algebras. This leads to formal moduli being given by differential
graded algebras, and gives rise the derived geometry programme. Usually, the
expectation is, that to solve a given 
{\em global }moduli problem with a differential graded Lie algebra, this
differential graded Lie algebra would have to infinite-dimensional,
and therefore ill-suited for algebraic geometry.  

For example, gauge theory can be used to construct analytic moduli
spaces of holomorphic vector bundles on a compact complex manifold $Y$.  In
this case, the differential graded Lie algebra might be
$A^{0,\bullet}(Y, M_n)$, the algebra of $C^\infty$-forms of type
$(0,\ast)$ with values in $n\times n$-matrices. The differential is
the Dolbeault differential, and the bracket is combined from wedge
product of forms and commutator bracket of matrices. Almost complex
structures are elements of $x\in A^{0,1}(Y,M_n)$, and they are integrable,
if and only if they satisfy the Maurer-Cartan equation
$$dx+\sfrac{1}{2}[x,x]=0\,.$$
Dividing the Maurer-Cartan locus by the gauge group
$G=A^{0,0}(Y,GL_n)$, we obtain the moduli space of holomorphic
bundles.

One central observation of this paper is that there exists a finite
dimensional analogue of this construction for moduli of coherent
sheaves on a smooth projective variety over $\cc$.  Derived moduli of
sheaves have been constructed before (see
\cite{DerQuot} or
\cite{ModuliDGCats}), but we believe it is a new observation that
there is a finite dimensional differential graded Lie algebra with an
algebraic gauge group, solving this moduli problem globally.  Simply
by virtue of being the space of Maurer-Cartan elements in a
differential graded Lie algebra up to gauge equivalence,  the
moduli space automatically comes with a derived, or differential graded,
structure.  

This construction also leads one immediately to the examination of Geometric
Invariant Theory stability for this algebraic gauge group action.
Thus, another
result of this paper it that GIT stability for our algebraic gauge
group action reproduces the standard notion of stability for sheaves.

Let
$Y$ be a smooth projective variety with homogeneous coordinate ring
$A$, and $\alpha(t)\in\qq[t]$ a numerical polynomial.

We present a construction of the derived moduli scheme
of stable sheaves  on $Y$, as a
Geometric Invariant Theory quotient of the
derived scheme of 
actions.  The derived scheme of actions, $\RAct$, was introduced by
Ciocan-Fontanine and Kapranov in \cite{DerQuot} as an auxiliary
tool in their construction of the derived scheme of quotients,
$\RQuot$. 

The basic idea is to describe a coherent sheaf $\fF$ on $Y$ with
Hilbert polynomial $\alpha(t)$ in terms
of the associated finite-dimensional graded $A$-module
$$\Gamma_{[p,q]}\fF=\bigoplus_{i=p}^q\Gamma\big(Y,\fF(i)\big)\,,$$
with dimension vector $\alpha|_{[p,q]}=\big(\alpha(p),\ldots,\alpha(q)\big)$, for $q\gg p\gg0$.
In fact,  for any open bounded family $\UU$ of sheaves with Hilbert
polynomial $\alpha(t)$ on $Y$, there exist $q\gg p\gg0$, such that
$$\Gamma_{[p,q]}:\UU\longrightarrow
\left(\parbox{10.5pc}{graded $A$-modules in $[p,q]$ with dimension vector
    $\alpha|_{[p,q]}$}\right)$$ 
is an open embedding of moduli functors (i.e., of stacks). 

We construct a finite-dimensional differential graded Lie algebra 
$$L=\bigoplus_{n=0}^{q-p}L^n$$
together with an algebraic gauge group $G$ (the Lie algebra of $G$ is $L^0$),
acting linearly on $L$, such that $MC(L)/G$, the quotient of the
solution set of the Maurer-Cartan equation
\begin{equation}\label{MC}
dx+\sfrac{1}{2}[x,x]=0\,,\qquad  x\in L^1
\end{equation}
by the gauge group, is equal to the set (or rather stack) of graded
$A$-modules concentrated in degrees between $p$ and $q$ with dimension
vector $\alpha|_{[p,q]}$, up to isomorphism. 

We do this by fixing a finite graded vector space $V$ of dimension
$\alpha|_{[p,q]}$. Then the degree 1 part of our differential graded
Lie algebra is essentially $L^1=\Hom_\gr(A,\End_\cc V)$, the space of degree
preserving $\cc$-linear maps from $A$ to $\End_\cc V$, and the
solutions to the Maurer-Cartan equation~(\ref{MC}) turn out to be
precisely the algebra maps $A\to\End_\cc V$, i.e. the structures of
graded $A$-module on $V$.  Taking the quotient by the gauge group $G=\GL_\gr(V)$ of graded
automorphisms of $V$ can be viewed as removing the choice of basis in $V$. 

Equivalently, a family of $A$-modules can be viewed as a graded vector bundle of
rank $\alpha|_{[p,q]}$, i.e., a $G$-torsor, endowed with an $A$-action. 
This approach to constructing (derived) moduli of $A$-modules in these
two steps, by first constructing moduli of  vector bundles, i.e., the 
stack $BG$,  and then a relative (derived) scheme of
actions over $BG$ is standard.  For example, To\"en and Vaqi\'e use this method
to construct moduli of derived category objects,
see~\cite{ModuliDGCats}. 

Our main interest lies in the derived {\em scheme}, obtained by
restricting to stable objects (which are simple) and then removing the
automorphism group (which is $\cc^\ast$), by passing to the space
underlying the $\cc^\ast$-gerbe.

The quotient $L^1/G$ is an instance of a moduli space of quiver
representations.  The relevant quiver is directed, which implies that
all points of $L^1$ are unstable for the action of $G$.  Using standard
techniques (as in \cite{QuiverKing}), we modify the action of $G$ on $L^1$
by a (canonical choice of) character of $G$ to obtain a well-defined
GIT (Geometric Invariant Theory) problem. Thus we obtain a quasi-projective moduli
space $MC(L)^s\git G$ of stable $A$-modules, with a compactification
$MC(L)^{ss}\git G$ consisting of semi-stable $A$-modules. The key result
is the 
following:

\paragraph{Theorem} For every  bounded family $\UU$ of sheaves on $Y$
with Hilbert polynomial $\alpha(t)$, there exist
$q\gg p\gg0$, such that if $\fF$ is a member of $\UU$, then $\fF$ is
a stable sheaf if and only if 
$\Gamma_{[p,q]}\fF$ is GIT-stable.

\bigskip

This shows that usual (semi)-stability as defined by
Simpson~\cite{simpson} is the natural notion of
(semi)-stability for sheaves induced from GIT-stability via our
construction.  Moreover, using the fact that semi-stable sheaves are
bounded, and satisfy the valuative criterion for properness, we see
that the moduli space of semi-stable sheaves with Hilbert polynomial
$\alpha(t)$ is a union of connected components of the projective
scheme $MC(L)^{ss}\git G$ of semi-stable modules.

This gives a new construction of the moduli space of (semi)-stable
sheaves on a projective variety. 
One advantage of our approach over others, such as the
classical Quot-scheme approach of \cite{simpson} and
\cite{HuybrechtsLehn},
or the Quiver approach of \cite{Alv-King}, is that (\ref{MC}) provides
us with a rather explicit set of equations cutting out the moduli
space.

We can also explicitly describe the image of the
moduli space of stable sheaves inside the moduli space of $[p,q]$-graded
$A$-modules.  Namely, it is the scheme of stable modules whose truncation into an interval
$[p',q]$, for suitable $p'$ between $p$ and $q$, is also stable.

Since $MC(L)/G$ is the moduli space of a differential graded Lie
algebra, it is automatically a differential graded
scheme. It is
naturally embedded into the smooth scheme $L^1/G$  as the `spectrum' of
a sheaf of differential graded algebras $\R$ on $L^1/G$,
obtained from the 
algebra of functions on the affine supermanifold 
$$ L[1]^{\geq0}\,,$$
with its induced derivation,
by descending to the $G$-quotient. It is this differential graded
scheme structure on $MC(L)/G$ which we refer to as a {\em derived scheme}. 

A derived scheme comes with higher obstruction spaces at every point.
In our case, the higher obstruction spaces at the sheaf $\fF$, or the
corresponding point $M=\Gamma_{[p,q]}\fF$ of $MC(L)/G$, are given by
$$\Ext^i_{\OO_Y}(\F,\F)=\Ext^i_A(M,M)\,.$$
The corresponding virtual fundamental class (see \cite{BF} and
\cite{vfcdgm}) is thus the one giving rise to Donaldson-Thomas
invariants~\cite{HolCas}. 

The differential graded Lie algebra $L$ is essentially the degree
preserving 
part of the Hochschild cochain complex 
$$L^n=\Hom_\cc(A^{\otimes n},\End_\cc V)$$
of the graded ring $A$ with values in the graded bimodule $\End_\cc V$, where $V$ is a
finite dimensional graded $A$-module in degrees from $p$ to $q$ with
dimension vector $\alpha|_{[p,a]}$, together with its natural Lie bracket
induced from the commutator bracket in $\End V$.

In Section~1, we construct the derived scheme of finite graded $A$-modules
with fixed dimension vector. This works for any algebra over $\cc$, in
particular, there is no need for commutativity of $A$.   The main purpose of
this section is to carefully describe the various differential graded
schemes and stacks we construct, and to do this as explicitly as
possible in terms of our finite
dimensional differential graded Lie algebra with its gauge group. We
hope the introduction of bundles of marked differential graded Lie
algebras will clarify the global geometric objects described
infinitesimally by differential graded Lie algebras. We also advocate
the use of Maurer-Cartan equations as convenient was to package higher
structures, in particular, $A_\infty$-module structures. 

Section~2 is devoted to the study of the GIT problem given by the
action of the gauge group $G$ on the space $L^1$. In particular, we
construct quasi-projective derived moduli spaces of equivalence
classes of stable finite graded $A$-modules of given dimension
vector. We hope there
will be applications in non-commutative geometry.

In Section~3 we introduce our projective scheme $Y$, and consider the
case where our graded ring $A$ is the homogeneous coordinate ring of
$Y$. We compare the stability notions for sheaves on $Y$ and for
graded $A$-modules.  We prove the above theorem and the amplification mentioned. 

Finally, in Section~4, we write down the derived moduli problem for
sheaves on $Y$, which
is solved by our differential graded scheme.  This is the only place
where we need $Y$ to be smooth.  The reason we need smoothness,
is to assure that for a coherent sheaf $\F$ on $Y$, the spaces
$\Ext^i_{\OO_Y}(\F,\F)$ vanish, for sufficiently large $i$. 

\subsubsection{Derived geometry}

For us, derived geometry is the geometry of differential graded
schemes. We make a few informal remarks here. For more detailed
expositions of derived geometry, see
To\"en-Vezzosi~\cite{BraveNewPaper}~\cite{HOGI} 
and~\cite{HOGII}  or Lurie~\cite{lurie}.

A {\em differential graded scheme }is a pair $(T,\R_T)$, where $T$ is a
scheme and $\R_T$ is a sheaf of differential graded $\cc$-algebras on
$T$, endowed with a structure morphism of sheaves of algebras
$\OO_T\to\R^0_T$.

It is natural to require (and we make it part of the definition), that
all differential graded schemes
$(T,\R_T)$ satisfy $\OO_T=\Z^0(\R_T)$, where $\Z^0(\R_T)=\ker(d:\R^0_T\to
  \R^1_T)$ is the sheaf of 0-cycles in $\R_T$. This implies that $\R_T$ is a sheaf of
differential graded $\OO_T$-algebras. Then  a {\em morphism }of
differential graded schemes $(T,\R_T)\to(M,\R_M)$ is a pair
$(\phi,\mu)$, where $\phi:T\to M$ is a morphism of schemes, and
$\mu:\phi^\ast \R_M\to \R_T$ is a morphism of sheaves of differential
graded $\OO_T$-algebras. 

The {\em classical scheme }associated to a differential graded scheme
$(T,\R_T)$ is the closed subscheme of $T$ given by
$\pi_0(T,\R_T)=\spec_{\OO_T}h^0(\R_T)$. 

A differential graded scheme is {\em affine}, if it comes from a
differential graded algebra which is free as a graded algebra, on a
finite set of generators, all in non-positive degree. 

Differential graded schemes form a category. (One may replace
morphisms by germs of morphisms, defined in suitable neighbourhoods of
the classical loci.) This category is enriched
over simplicial sets: the $n$-simplices in $\Hom(X,Y)$ are
the morphisms $X\times\Delta_n\to Y$, where $\Delta_n$ is the
differential graded scheme (which is not affine) corresponding to the differential graded
algebra of algebraic differential forms on the algebraic
$n$-simplex.

The category of differential graded schemes also has a
natural topology: the \'etale topology, in which a family $U_i\to U$
is a covering family if $\pi_0(U_i)\to \pi_0(U)$ is a covering family
in the usual \'etale topology, and every $U_i\to U$ is an \'etale
morphism, which means that
$h^r(\R_{U_i})=h^r(\R_{U})\otimes_{h^0(\R_{U})}h^0(\R_{U_i})$, for all $r$. 

In analogy with the definition of algebraic spaces, one can define a
derived scheme (or space) to be a simplicial presheaf $X$ on
the category of differential graded schemes satisfying two properties:
\begin{items}
\item (sheaf property) for every hypercover $U_\bullet\to U$ the map
  $X(U)\to \mathop{\rm hocolim} X(U_\bullet)$ is a weak equivalence,
\item (locally affine property) $X$ is \'etale locally weakly
  equivalent to a presheaf represented by an affine differential
  graded scheme. 
\end{items}

A particularly nice class of differential graded schemes comes from
bundles of marked differential graded Lie algebras on smooth schemes 
(see the beginning of Section~\ref{12}). Our main object of study,
$\RModtilde_{\alpha}^\sp(A)=(\tilde\MM^\sp,\R)$ is of this form.

We find it plausible (this will be proved elsewhere), that
differential graded schemes coming from bundles of marked differential
graded Lie algebras represent simplicial presheaves satisfying the
above two properties, (at least if we restrict the underlying base
category to affine objects).  Therefore, The moduli functor
represented by such a `nice' differential graded scheme, in the derived
world, would be given directly by the functor it represents over the
category of differential graded schemes as defined here.  This is the
moduli functor we examine.

\subsubsection{Glossary}

\begin{description}[topsep=3pt, itemsep=3pt, parsep=0pt,
  labelwidth=5pc, labelindent=5pc, leftmargin=5pc]
\item[$A$] A graded ring.
\item[$\m$] The maximal ideal of positive degree elements in $A$.
\item[$\alpha$] Depending on the context, either a numerical
  polynomial $\alpha(t)\in\qq[t]$, or a dimension vector
  $(\alpha_p,\ldots,\alpha_q)$. 
\item[$V$] A graded vector space of dimension $\alpha=(\alpha_p,\ldots,\alpha_q)$.
\item[$L$] The differential graded Lie algebra $L^n=\Hom_\gr(\m^{\otimes
    n},\End V)$, see Section~\ref{dgll}.
\item[$M$] The scheme $L^1$. 
\item[$X$] The Maurer-Cartan locus in $M$. 
\item[$\R_M$] The sheaf of differential graded algebras on $M$, see
  Section~\ref{12}.
\item[$\Act_\gr(A,V)$] The scheme $X$, when it is viewed as
  representing the scheme of graded actions of $A$ on $V$.
\item[$\RAct_\gr(A,V)$] The differential graded scheme $(M,\R_M)$,
  which is the derived scheme of actions. 
\item[$G$] The gauge group $G=\prod_{i=p}^q GL(V_i)$. 
\item[$\Delta$] The one-parameter subgroup of scalars in $G$. 
\item[$\tilde G$] The quotient group $G/\Delta$. 
\item[$\MM$] The quotient stack $[M/G]$
\item[$\tilde\MM$] The quotient stack $[M/\tilde G]$. 
\item[$\tilde\MM^\sp$] The open substack of $\tilde\MM$ which is an
  algebraic space.
\item[$\XX$] The Maurer-Cartan locus in $\MM$.
\item[$\tilde\XX$] The Maurer-Cartan locus in $\tilde\MM$. 
\item[$\Mod_\alpha(A)$] The algebraic stack $\XX$, when it is viewed
  as the stack of graded $A$-modules of dimension $\alpha$. 
\item[$\RMod_\alpha(A)$] The differential graded stack
  $(\MM,\R_{\MM})$, which is the derived stack of graded modules.
\item[$\Modtilde_\alpha^\sp(A)$] The algebraic space $\tilde\XX$, when
  viewed as the space of equivalence classes of simple graded modules.
\item[$\RModtilde_\alpha^\sp(A)$] The differential graded algebraic
  space $(\tilde\MM^\sp,\R)$, which is the derived space of
  equivalence classes of simple modules. 
\item[$\RModtilde_\alpha^s(A)$] The stable locus inside
  $\RModtilde_\alpha^\sp(A)$.
\item[$\RModtilde_{\alpha}^\sp(\OO_Y)$] The functor of
  equivalence classes of simple families of coherent sheaves on $Y$
  with Hilbert polynomial $\alpha(t)$ parametrized by differential graded schemes.
\item[$\RModtilde_{\alpha}^s(\OO_Y)$] The stable locus inside
  $\RModtilde_{\alpha}^\sp(\OO_Y)$.
\end{description}

\subsubsection{Notation and Conventions}

We work over a field of characteristic zero, which we shall denote by
$\cc$.  All tensor products are over $\cc$, unless indicated
otherwise.  All our differential graded algebras (and sheaves
thereof), are graded commutative with unit. 

Cohomology sheaves (of a complex of sheaves $\E^\com$) we
usually denote by~$h^i(\E)$.

\section{The derived scheme of simple graded modules} \label{RAct}

Let $A$ be a unital graded $\cc$-algebra, {\em not necessarily
  commutative}, which is all in non-negative degrees, and such that each
graded piece is finite dimensional.  Moreover, we assume that the
degree zero piece is one-dimensional, hence equal to $\cc$.  We
denote by $\m$ the ideal of elements of positive degree in $A$.  Note that $\m$
is a positively graded algebra without unit.  We refer to the grading on $A$ as
the {\em internal }or {\em projective }grading, if there is a fear of
confusion.  We indicate this grading with lower indices.

Our main example of interest is that $A$ is the
homogeneous coordinate ring of a projective variety over $\cc$.

A graded $A$-module is the same thing as a graded $\m$-module.  The
advantage of working with $\m$ is that there is only one module axiom:
associativity.

\subsection{The \dg Lie algebra $L$}\label{dgll}

Let $V$ be a graded and finite-dimensional vector space
$$V=\bigoplus_{i=p}^q V_i\,.$$
By $\End V$ we denote the algebra of $\cc$-linear endomorphisms of
$V$. It inherits a grading from $V$. Only $\End_iV$ in the range
$i\in[p-q,q-p]$ are non-zero.

By
$$\alpha=(\alpha_p,\ldots,\alpha_q)=(\dim V_p,\ldots,\dim V_q)$$
we denote the dimension vector of $V$.

\subsubsection{The graded vector space}

We consider
$$L^n = \Hom_\gr(\m^{\otimes n}, \End V)\,,$$
the vector space of degree-preserving $\cc$-linear maps
$\mu:\m^{\otimes n}\to\End V$, and
$$L=\bigoplus_{n=0}^\infty L^n\,.$$
Thus, $L^0=\End_\gr V$ and $L^1=\Hom_\gr(\m,\End V)$.
We write elements $\mu\in L^n$ as multilinear maps $\m^{\times n}\to
\End V$. To distinguish the grading on $L$ from the projective
grading, we may sometimes refer to it as the {\em external
  grading}. It is always indicated by upper indices.

Note that every $L^n$ is finite-dimensional, and that $L^n=0$, unless
$n$ is in the range $n\in[0,q-p]$, because $\m$ is positively graded.

Each $L^n$ is graded projectively:
$$L^n=\bigoplus_{\substack{q\geq i\geq j\geq p\\[2pt] i-j\geq n}}
L^n_{ij}\,,$$
where
$$L_{ij}^n=
\Hom\big(
  (\m^{\otimes n})_{i-j}, \Hom(V_j,V_i)\big)\,.$$
For $n=0$, this simplifies to
$$L^0=\bigoplus_{i=p}^q L^0_{ii}\,,\qquad L^0_{ii}=\Hom(V_i,V_i)\,,$$
and for $n=1$, we can write
$$L^1=\bigoplus_{q\geq i>j\geq p} L^1_{ij}\,,\qquad
L^1_{ij}=\Hom\big( \m_{i-j} , \Hom(V_j,V_i)\big)\,.$$
We say that $L^0$ is diagonal, and $L^1$ is strictly lower
triangular. The higher $L^n$ are restricted to successively smaller
`South-West' corners.

\subsubsection{The gauge group}

We let $G=GL_\gr(V)$ be the group of degree-preserving linear
automorphisms of $V$, and call it the {\em gauge group}. Of course,
$L^0$ is the Lie algebra of $G$. The gauge group is graded:
$$G=\prod_{i=p}^q G_i\,,\qquad  G_i=GL(V_i)\,.$$
It acts, from the left,  via conjugation on $L$. More precisely, for $g\in G$ and
$\mu\in L^n$, we have
\begin{equation}\Label{act}
(g\cdot\mu)(a_1,\ldots,a_n)=g\circ\mu(a_1,\ldots,a_n)\circ
g^{-1}\,.
\end{equation}
The action of $G$ on $L^n$ preserves the double grading: if
$g=(g_p,\ldots,g_q)$, and $\mu\in L^n$, then
\begin{equation}\Label{actij}
(g\cdot\mu)_{ij}=g_i\mu_{ij}g_j^{-1}\,.
\end{equation}
We call this action the {\em gauge action}. 
The group $G$ contains the scalars, $\Delta:\cc^\ast\to G$,
$t\mapsto(t,\ldots,t)$, which act trivially.  This leads
us to also consider the quotient group $\tilde G=G/\Delta$. 

\subsubsection{The differential}

Define $d:L^n\to L^{n+1}$ by the formula
$$d\mu(a_1,\ldots,a_{n+1})=\sum_{i=1}^n(-1)^{n-i}\mu(\ldots,a_ia_{i+1},\ldots)$$
For example, $d:L^0\to L^1$ is equal to zero, and $d:L^1\to L^2$ is
given by $d\mu(a,b)=\mu(ab)$.

Of course, $d^2=0$. The gauge action preserves the differential. The
differential preserves the projective double grading.

The complex $(L,d)$ is the subcomplex of internal degree zero of the
Hochschild complex of 
the $\cc$-algebra $\m$ with values in the bimodule $\End V$, where
$\End V$ has the trivial (i.e., zero) module structure.

\subsubsection{The bracket}

For $\mu\in L^m$ and $\mu'\in L^n$ define $\mu\circ\mu'\in L^{m+n}$ by the formula
$$\mu\circ\mu'\,(a_1,\ldots,a_{m+n})=
(-1)^{mn}\mu(a_1,\ldots,a_m)\circ\mu'(a_{m+1},\ldots,a_{m+n})\,.$$ 
An easy sign calculation shows that this operation is associative.

Then, for $\mu\in L^m$ and $\mu'\in L^n$ define $[\mu,\mu']\in
L^{m+n}$ by
$$[\mu,\mu']=\mu\circ\mu'-(-1)^{mn}\mu'\circ \mu\,.$$
This operation { automatically }satisfies the graded Jacobi identity,
because it is defined as the graded commutator of an associative
product.

We can write out the formula for the bracket:
\begin{multline*}
[\mu,\mu'](a_1,\ldots,a_{m+n})=
(-1)^{mn}\mu(a_1,\ldots,a_m)\circ\mu'(a_{m+1},\ldots,a_{m+n}) \\
-\mu'(a_1,\ldots,a_n)\circ\mu(a_{n+1},\ldots,a_{m+n})\,.
\end{multline*}
For example, if $\mu,\mu'\in L^1$, then
$$[\mu,\mu'](a,b)=-\mu(a)\circ\mu'(b)-\mu'(a)\circ\mu(b)\,.$$

The differential $d$ acts as a derivation with respect to the
bracket~$[\,]$, i.e., for $\mu\in L^m$ and $\mu'\in L^n$, we have
$$d[\mu,\mu']=[d\mu,\mu']+(-1)^m[\mu,d\mu']\,.$$
Thus $(L,d,[\,])$ is a differential graded Lie algebra.

The gauge group $G$ acts by automorphisms of the \dg Lie algebra structure on
$L$.
This means that we have
$$d(g\cdot\mu)=g\cdot d\mu\qquad\text{and}\qquad
g\cdot[\mu,\mu']=[g\cdot\mu,g\cdot\mu']\,.$$ 

The derivative of the gauge action of $G$ on $L$ is the adjoint action
of $L^0$ on $L$.

\begin{rmk}
The more basic object than $L$ is the truncation $L^{>0}=\tau_{>0}L$, together
with $G$ and its gauge action. The differential graded Lie algebra $L$ can
be recovered from $(L^{>0},G)$.
\end{rmk}

\subsubsection{The Maurer-Cartan equation}

The Maurer-Cartan equation is
$$d\mu+\tfrac{1}{2}[\mu,\mu]=0\qquad\text{for $\mu\in L^1$}\,.$$
We call $\mu\in L^1$ a {\em Maurer-Cartan element}, if it satisfies this
equation. We denote the set of Maurer-Cartan elements by $MC(L)$.

For $\mu\in L^1$, we have
$\tfrac{1}{2}[\mu,\mu]=\mu\circ\mu$,
and so $\mu$ is a Maurer-Cartan element if and only if
$$d\mu+\mu\circ\mu=0\,,$$
or, equivalently, if for
all $a,b\in \m$,
$$\mu(ab)=\mu(a)\circ\mu(b)\,.$$
If we write out this equation degree-wise, we
get  for all $i>k>j$ and $a\in \m_{i-k}$, $b\in \m_{k-j}$ 
the equation
$\mu_{ij}(ab)=\mu_{ik}(a)\circ\mu_{kj}(b)$.

Thus $\mu\in L^1$ is a Maurer-Cartan element if and only if it defines
a { left }action of $\m$ on $V$.  Dividing by the gauge action removes
the choice of basis in $V$. It follows immediately, that Maurer-Cartan
elements up to gauge equivalence are graded $\m$-modules up to isomorphism, whose
underlying graded vector space is isomorphic to $V$. We can make this
claim precise:

\begin{rmk}
Let $[MC(L)/G]$ be the (set-theoretic) transformation groupoid associated
to the gauge group action on the Maurer-Cartan elements. Let
$(\text{$\m$-modules})_{\alpha}$ denote the category of graded
$\m$-modules with dimension vector $\alpha$, with only
isomorphisms. Then we have an equivalence of groupoids
$$[MC(L)/G]\longrightarrow(\text{$\m$-modules})_{\alpha}\,,$$
given by mapping $\mu$ to the $\m$-module structure it defines on $V$
and mapping an element of $G$ to the isomorphism of $\m$-module
structures it represents. We will turn this into a
geometric statement. 
\end{rmk}

\subsection{The moduli stack of $L$}\label{12}

The following construction of the differential graded moduli stack
works for any 
finite-dimensional differential graded Lie algebra  concentrated in
non-negative degrees with algebraic gauge group.

\subsubsection{Bundles of marked differential graded Lie algebras}

\begin{defn}\label{bom}
A {\em bundle of marked differential graded Lie algebras }over a
scheme (or a stack) $M$ is a graded vector bundle $\LL^\ast$ over $M$, endowed with
three pieces of data: 
\begin{items}
\item a section $f\in\Gamma(M,\LL^2)$,
\item an $\OO_M$-linear map of degree one $\delta:\LL^\ast\to\LL^\ast$,
\item a $\OO_M$-linear alternating bracket of degree zero
  $[\,]:\Lambda^2\LL^\ast\to\LL^\ast$,
\end{items}
subject to four axioms:
\begin{items}
\item $\delta(f)=0$, as a section of $\LL^3$,
\item $\delta\circ \delta=[f,\,]$,
\item $\delta$ is a graded derivation with respect to the bracket $[\,]$,
\item the bracket $[\,]$ satisfies the graded Jacobi identity.
\end{items}
\end{defn}

A bundle of marked differential graded Lie algebras is a bundle of
differential graded Lie algebras only if $f=0$. All of our bundles of
marked differential graded Lie algebras will be concentrated in
degrees $\geq2$. The map $\delta$ will be referred to as the {\em
  twisted differential}. 

It will be useful to relax the conditions somewhat, and call a {\em
  sheaf }of marked differential graded Lie algebras on $M$ a graded
sheaf of $\OO_M$-algebras, endowed with the same data (i) to (iii),
subject to the same constraints (i) to (iv).  Sheaves of marked
differential graded Lie algebras will also be allowed to have
contributions in degrees less than 2. The sheaf of Maurer-Cartan
elements of a sheaf of marked differential graded Lie algebras is the
preimage of $-f$ under the curvature map $\LL^1\to\LL^2$ given by
$x\mapsto \delta x+\frac{1}{2}[x,x]$.  If $\LL$ is a
bundle (so that $\LL^1=0$), then the Maurer-Cartan locus is the
scheme-theoretic vanishing locus of $f$ in $M$.

If $\LL$ is a bundle of marked differential graded Lie algebras on $M$ and
$\R_M$ a sheaf of differential graded $\OO_M$-algebras, then 
$\LL\otimes_{\OO_M}\R_M$ is in a natural way a sheaf of marked
differential graded Lie algebras. 

\subsubsection{Associated differential graded scheme or stack}

To a bundle of marked differential graded Lie algebras over $M$ we associate a
sheaf of differential graded  algebras $\R_M$ by letting
the underlying sheaf of graded $\OO_M$-algebras be
\begin{equation}\label{sfl}
\R^\ast_M=\Sym_{\OO_M}\LL[1]^\vee\,,
\end{equation}
the sheaf of free graded commutative $\OO_M$-algebras with unit on the shifted
dual of $\LL$.

The bracket defines a morphism
$q_2:\LL[1]^\vee\to\Sym^2_{\OO_M}\LL[1]^\vee$, the twisted differential a
morphism $q_1:\LL[1]^\vee\to\Sym^1_{\OO_M}\LL[1]^\vee=\LL[1]^\vee$ and the
marking a morphism
$q_0:\LL[1]^\vee\to\Sym^0_{\OO_M}\LL[1]^\vee=\OO_M$. All three morphisms
$q_i$ have homological degree $+1$, and all three extend uniquely to
$\OO_M$-linear derivations $q_i:\R_M\to\R_M$.  Let $q=q_0+q_1+q_2$ be
the sum of these three derivations. The four axioms of marked
differential graded Lie algebra translate into the one condition 
$$q^2=0$$
for the derivation $q$ on $\R_M$.  This defines the differential
graded scheme $(\R_M,q)$. We will usually suppress $q$ from the notation.

Note that $X=Z(f)\subset M$, the scheme theoretic vanishing locus of
$f$ (the Maurer-Cartan locus), is equal to the subscheme of $M$
defined by the image of $\R^{-1}$ in $\R^0_M=\OO_M$. The structure
sheaf of $X$ is $\OO_X=h^0(\R_M)$.

\begin{ex}\label{lie}
Given a finite dimensional differential graded Lie algebra $L$,
concentrated in degrees $>0$, we let $M=L^1=\spec\Sym(L^{1\vee})$. Over
$M$ we consider for every $i\geq2$ the trivial vector bundle $\LL^i$
with fibre $L^i$, i.e., $\LL^i=L^i\times M$. The curvature map
$f:L^1\to L^2$ given by $f(x)=dx+\frac{1}{2}[x,x]$ gives rise to a
section of $\LL^2$ over $M$, the twisted differential $\delta=d^\mu:\LL^i\to
\LL^{i+1}$ is defined by the formula $\delta(y)=d^\mu(y)=dy+[\mu,y]$ in the fiber
over $\mu\in M=L^1$, and the bracket on $\LL$ is constant, i.e., equal to
the bracket on $L$ in every fibre of $\LL$. In this way the
differential graded Lie algebra $L=L^{\geq1}$ gives rise to a bundle of
differential graded Lie algebras $\LL=\LL^{\geq2}$ over $M=L^1$. 

Note that $X=Z(f)\subset M$ is identified with the scheme theoretic
Maurer-Cartan locus of $L$. 

If an algebraic group $G$ acts on $L$ by automorphisms of the
differential graded Lie algebra structure, the bundle of marked
differential graded Lie algebras $\LL$ over $M$ inherits a $G$ action
covering the $G$-action on $M$. Thus, the bundle of marked
differential graded Lie algebras $\LL$ descends to the quotient stack
$[M/G]$. 
\end{ex}

We apply these considerations to the truncation of our differential
graded Lie algebra $L^{>0}$ with the gauge group action by $G$.  We
obtain a bundle of marked differential graded Lie algebras $\LL_{\MM}$ over
$\MM=[M/G]$ and a sheaf of differential graded algebras $\R_{\MM}$ over
$\MM$. 

If we replace $G$ by $\tilde G$, we obtain a bundle of marked
differential graded Lie algebras $\LL_{\tilde \MM}$ over
$\tilde\MM=[M/\tilde G]$, and a sheaf of differential graded algebras
$\R_{\tilde\MM}$ over $\tilde\MM$. 
The Maurer-Cartan locus $X=Z(f)\subset M$ descends to closed substacks
$\XX\subset\MM$ and $\tilde\XX\subset\tilde\MM$, such that
$\OO_\XX=h^0(\R_\MM)$ and $\OO_{\tilde\XX}=h^0(\R_{\tilde\MM})$.

\begin{rmk}\label{13}
There is a natural morphism $\MM\to\tilde\MM$, making $\MM$ a
$\cc^\ast$-gerbe over $\tilde\MM$. This gerbe is trivial, if there
exists a line bundle $\xi$ over $M$ and a lifting of the $G$-action to
a $G$-action on $\xi$, such that $\Delta$ acts by scalar
multiplication on the fibres of $\xi$. 
\end{rmk}

\subsubsection{The associated functor on dg schemes}

Suppose the differential graded scheme $(M,\R_M)$ comes from a bundle of
marked differential graded Lie algebras as in (\ref{sfl}). Given a
morphism of schemes $\phi:T\to M$, the sheaf of Maurer-Cartan elements
of $\phi^\ast\LL\otimes_{\OO_T}\R_T$ is naturally isomorphic to the sheaf of
morphisms of differential graded $\OO_T$-algebras
$\phi^\ast\R_M\to\R_T$. 
$$MC(\phi^\ast\LL\otimes_{\OO_T}\R_T)=\Mor_{\OO_T}(\phi^\ast\Sym_{\OO_M}\LL[1]^\vee,\R_T)\,.$$
In particular, a 
morphism of differential graded schemes $(T,\R_T)\to(M,\R_M)$ is
essentially the same thing as a pair $(\phi,\mu)$, where $\phi:T\to M$
is a morphism of schemes and $\mu$ is a global Maurer-Cartan element
of the sheaf of marked differential graded Lie algebras
$\phi^\ast\LL\otimes_{\OO_T}\R_T$. 

\begin{lem}
If $(M,\R_M)$ comes as in Example~\ref{lie} from a differential graded
Lie algebra $L=L^{\geq 1}$, then a morphism $(T,\R_T)\to(M,\R_M)$ is the same
thing as a 
global Maurer-Cartan element in the sheaf of differential graded Lie
algebras $L\otimes_\cc\R_T$. 
\end{lem}
\begin{pf}
A morphism $\phi:T\to L^1$ can be considered as a section of
$L^1\otimes\Z^0(\R_T)$ and hence as a degree $1$ section of
$L\otimes\R_T$. The section $\mu$ can also be thought of as a degree
$1$ section $L\otimes \R_T$, and it is not hard to check, that $\mu+\phi$ is a
Maurer-Cartan section of the sheaf of 
differential graded Lie algebras $L\otimes\R_T$. Conversely, every
Maurer-Cartan section of $L\otimes \R_T$ gives rise to a pair
$(\phi,\mu)$ and hence to a morphism of differential graded schemes
$(T,\R_T)\to(M,\R_M)$. 
\end{pf}

Finally, if $G$ acts on $L$ by automorphisms, and $\MM=[M/G]$, then a
morphism $(T,\R_T)\to(\MM,\R_\MM)$ is {\em essentially }the same thing as a pair
$(E,\mu)$, where $E$ is a principal $G$-bundle over $T$, and $\mu$ is a
global Maurer-Cartan element of the sheaf of differential graded Lie
algebras ${}^EL\otimes_{\OO_T}\R_T$. Here ${}^EL$ denotes the
associated vector bundle, with its induced structure of sheaf of
differential graded Lie algebras over $\OO_T$.

We apply these considerations to the differential graded Lie algebra
$$L^{\geq1}=\Hom_\gr(\m^{\otimes\geq1},\End V)\,.$$

\subsubsection{The derived scheme of actions}

Let $(M,\R_M)$ be the differential graded scheme associated as in
Example~\ref{lie} to $L^{\geq1}=\Hom_\gr(\m^{\otimes\geq1},\End V)$.  So
$M=\Hom_\gr(\m,\End V)$. The following is essentially
Proposition~(3.5.2) of \cite{DerQuot}.

\begin{prop}
Let $(T,\R_T)$ be a differential graded scheme. A morphism
$(T,\R_T)\to(M,\R_M)$ is essentially the same thing as a Maurer-Cartan
element in the differential graded Lie algebra
$$\Gamma\big(T,\Hom(\m^{\otimes\geq1},\End V)\otimes\R_T\big)\,.$$
This, in turn, is the same thing as a graded $\R_T$-linear
$A_\infty$-action of $\m\otimes\R_T$ on $V\otimes\R_T$, or a graded
{\em unitary }$\R_T$-linear $A_\infty$-action of $A\otimes\R_T$ on $V\otimes \R_T$. 
\end{prop}

This justifies calling $(M,\R_M)$ the {\em derived scheme of graded
  actions }of $A$ on $V$, and denoting it by $\RAct_\gr(A,V)$. 

\subsubsection{The derived stack of modules}

Let $(\MM,\R_{\MM})$ be the differential graded stack obtained from
$(M,\R_M)$ by dividing by $G$, and let $\XX\subset\MM$ be the
Maurer-Cartan locus. 

\begin{prop}
Let $(T,\R_T)$ be a differential graded scheme. A morphism
$(T,\R_T)\to(\MM,\R_\MM)$ is essentially the same thing as a pair
$(E,\mu)$, where $E=\bigoplus_{i=p}^q E_i$
is a graded vector bundle of dimension vector $\alpha$ over $T$, and
$\mu$ is a Maurer-Cartan element in the differential graded Lie
algebra
$$\Gamma\big(T,\sheafhom_\gr(\m^{\otimes\geq 1},\sheafend_{\OO_T}
E)\otimes_{\OO_T}\R_T\big)\,.$$
Such a Maurer-Cartan element $\mu$ is the same thing as a graded
$\R_T$-linear $A_\infty$-action of $\m\otimes\R_T$ on
$E\otimes_{\OO_T}\R_T$, or a graded {\em unitary }$\R_T$-linear
$A_\infty$-action of $A\otimes\R_T$ on $E\otimes_{\OO_T}\R_T$. 

In particular, if $T$ is a classical scheme, a morphism
$T\to(\MM,\R_\MM)$ is the same thing as a morphism $T\to\XX$, which,
in turn, is the same thing as a graded vector bundle over $T$ of
dimension $\alpha$, endowed with the structure of a sheaf of graded
$\m\otimes\OO_T$-modules, or a the structure of a sheaf of graded {\em
  unitary }$A\otimes\OO_T$-modules.
\end{prop}

There is a universal family over $(\MM,\R_{\MM})$.  It is obtained
from $V\otimes\OO_M$ with its tautological $A_\infty$-action
$$\boldsymbol{\mu}:\m\otimes V\otimes\R_M\longrightarrow V\otimes\R_M\,,$$
by descent: the group $G$ acts naturally on $V$, in a way respecting
$\boldsymbol{\mu}$. 

We call $(\MM,\R_\MM)$ the {\em derived stack of graded $A$-modules }with
dimension vector $\alpha$, and use the notation
$\RMod_\alpha(A)=(\MM,\R_{\MM})$. For the underlying classical stack
$\XX$ we write $\Mod_\alpha(A)=\XX$.

\subsection{The derived space of equivalence classes of simple modules}

When dividing by $\tilde G$ instead of $G$, we have to be more
careful, because the natural action of $G$ on $V$ does not factor
through $\tilde G$, as the scalars in $G$ do not act trivially on $V$.
This implies that the universal family
of graded $A$-modules does not descend from $\MM$ to $\tilde\MM$. 
The obstruction is the $\cc^\ast$-gerbe of Remark~\ref{13}.

\subsubsection{Equivalence of simple modules}

A {\em family of graded $A$-modules }of dimension $\alpha$
parametrized by the scheme $T$, is a graded vector bundle with rank
vector $\alpha$ on $T$ together with a unitary graded $\OO_T$-linear action of
$A\otimes\OO_T$. 

\begin{defn}
A family $E$ of graded $A$-modules parametrized by $T$ is {\em
  simple}, if the sheaf of endomorphisms of $E$ is equal to $\OO_T$. 
Two simple families of graded $A$-modules $E$, $F$, parametrized by $T$ are {\em
  equivalent}, if there exists a line bundle $\lL$ on $T$, such that
$F$ is isomorphic to $E\otimes_{\OO_T}\lL$, as family of graded
$A$-modules.
\end{defn}

Equivalence classes of simple families of graded $A$-modules form a presheaf
on the site of $\cc$-schemes with the \'etale topology, whose
associated sheaf we denote by $\Modtilde_\alpha^\sp(A)$. 

Let $M^\sp\subset M$ be the open subscheme of points with trivial $\tilde
G$-stabilizer, and $X^\sp=X\cap M$ the intersection with the
Maurer-Cartan locus $X$.  Denote
by $\tilde\MM^\sp\subset\tilde\MM$ and $\tilde\XX^\sp\subset\tilde
\XX$ the quotients by $\tilde G$. 

\begin{rmk}
The sheaf $\Modtilde_\alpha^\sp(A)$ is isomorphic to the algebraic
space~$\tilde\XX^\sp$.  
$$\Modtilde_\alpha^\sp(A)=\tilde\XX^\sp$$
This proves that $\Modtilde_\alpha^\sp(A)$ is
algebraic, and gives a modular interpretation of $\tilde\XX^\sp$. 
\end{rmk}

\subsubsection{Coprime case}

\begin{prop}
Suppose that the components of the dimension vector $\alpha$ are
coprime.  Then the gerbe of Remark~\ref{13} is trivial. Moreover, the
presheaf of equivalence classes of simple families 
of graded $A$-modules is a sheaf.  In other words, for any
$\cc$-scheme $T$, the $T$-points of the algebraic space
$\Modtilde_\alpha^\sp(A)$ correspond one-to-one to equivalence
classes of simple families. In particular,
$\Modtilde_\alpha^\sp(A)$ admits a universal family of simple
graded $A$-modules. 
\end{prop}
\begin{pf}
There exist integers $n_i$ such that
$\sum_{i=p}^qn_i\alpha_i=1$.
The character $\rho:G\to\cc^\ast$ given by 
$\rho(g)=\prod_{i=p}^q\det(g_i)^{n_i}$
satisfies $\langle\Delta,\rho\rangle=1$. So twisting the action of $G$ on
$V$ by $\rho^{-1}$, the twisted action factors
through $\tilde G$, and so after the twist, $V$ descends to
$\tilde\MM$. 
\end{pf}

\begin{rmk}
If 
$$\alpha(t)=a_0\binom{t}{0}+a_1\binom{t}{1}+\ldots+a_k\binom{t}{k}$$
is a numerical polynomial $\alpha(t)\in\qq[t]$ of degree $k$, with
$a_0,\ldots,a_k\in\zz$, 
and $q-p\geq k$, then
$$\big(\alpha(p),\ldots,\alpha(q)\big)=1\Longleftrightarrow(a_0,\ldots,a_k)=1\,.$$
Hence $\big(\alpha(p),\ldots,\alpha(q)\big)=1$ if and only if $\alpha$
is primitive (not an integer multiple of another numerical
polynomial). 
\end{rmk}

We will write down the derived moduli problem solved by the
differential graded algebraic space $(\tilde\MM^\sp,\R)$.

\subsubsection{The derived space of simple modules}

Let $(T,\R_T)$ be a differential graded scheme. If $F$ is a 
graded vector bundle on $T$, we can sheafify the
construction of our differential graded Lie algebra over $T$, and
tensor with $\R_T$  to obtain
a sheaf of differential graded Lie algebras 
\begin{equation}\label{djksfl}
\sheafhom_\gr(\m^{\otimes\geq1},\sheafend_{\OO_T}F)\otimes_{\OO_T}\R_T\,.
\end{equation}
A global Maurer-Cartan element in (\ref{djksfl}) is the same thing as
a graded
$\R_T$-linear $A_\infty$-action of $\m\otimes\R_T$ on
$F\otimes_{\OO_T}\R_T$.

A {\em family of graded $A$-modules }with dimension vector $\alpha$
parametrized by the differential graded scheme $(T,\R_T)$ is a pair
$(F,\mu)$, where $F$ is a graded vector bundle of dimension $\alpha$
over $T$, and $\mu$ is a global Maurer-Cartan element
in~(\ref{djksfl}). Two such families are {\em equivalent}, if they
differ by a line bundle on $T$. We denote the set of equivalence
classes of such families by $\RModtilde_\alpha(A)(T)$. Varying
$(T,\R_T)$, we get a presheaf $\RModtilde_\alpha(A)$ on the category of
differential graded schemes.

Note that a Maurer-Cartan element $\mu$ in (\ref{djksfl}) can be
decomposed
$$\mu=\sum_{i=1}^{q-p} \mu_i\,,\qquad \mu_i\in\sheafhom_\gr(\m^{\otimes
  i},\sheafend_{\OO_T}F)\otimes_{\OO_T}\R^{1-i}_T\,.$$ So
$\mu_1\in\sheafhom_\gr(\m,\sheafend_{\OO_T}F)\otimes_{\OO_T}\R^{0}_T$. The
Maurer-Cartan equation implies that $\mu_1$ takes values in the
subsheaf
$\sheafhom_\gr(\m,\sheafend_{\OO_T}F)\otimes_{\OO_T}\Z^0(\R_T)$, which
is equal to $\sheafhom_\gr(\m,\sheafend_{\OO_T}F)$, by our definition
of differential graded scheme. Thus, we may also think of $\mu_1$ as an
$\OO_T$-linear map
$\mu_1:\m\otimes\OO_T\to\sheafend_{\OO_T}F$.
We call $(F,\mu)$ {\em simple}, if the subsheaf of
$\sheafend_{\OO_T}F$ commuting with the image
of $\mu_1$ is equal to $\OO_T$.  Simple families define the
subpresheaf $\RModtilde_\alpha^\sp(A)\subset\RModtilde_\alpha(A)$.

\begin{prop}
The differential graded algebraic space $(\tilde\MM^\sp,\R)$
represents the sheaf associated to $\RModtilde_\alpha^\sp(A)$. If
$\alpha$ is primitive, then $\RModtilde_\alpha^\sp(A)$ is a sheaf, and
so $(\tilde\MM^\sp,\R)$ represents $\RModtilde_\alpha^\sp(A)$. 
\end{prop}
\begin{pf}
Let $(F,\mu)$ be a simple graded family parametrized by the
differential graded scheme $(T,\R_T)$. Write $\mu=\mu_1+\mu'$, where
$\mu'=\sum_{i\geq2}\mu_i$. Then the pair $(F,\mu_1)$ defines a
morphism $\phi:T\to\tilde\MM^\sp$, and any equivalent simple graded family
gives rise to the same morphism $T\to\tilde\MM^\sp$. 
The pullback to $T$
of $\R_{\tilde \MM}$ via the morphism $\phi$ is equal to the sheaf of
symmetric algebras generated over $\OO_T$ by the shifted dual of 
$\sheafhom_\gr(\m^{\geq2},\sheafend_{\OO_T} F)$. Therefore, a morphism
$\phi^\ast\R_{\tilde\MM}\to\R_T$ is the same thing as a global Maurer-Cartan
section of the sheaf of marked differential graded Lie algebras
(with twisted differential)
$$\sheafhom_\gr(\m^{\geq2},\sheafend_{\OO_T} F)\otimes_{\OO_T}\R_T\,.$$
This is exactly what $\mu'$ provides us with. Hence $(F,\mu)$ gives
rise to a morphism $(T,\R_T)\to(\tilde\MM,\R)$. 

We have defined a morphism from the presheaf 
$\Modtilde_\alpha^\sp(A)$
to the sheaf represented by $(\tilde\MM^\sp,\R)$. 
Conversely, every morphism $\phi:T\to\tilde\MM^\sp$ is (locally in $T$),
induced by a pair $(F,\mu_1)$, and every morphism
$\phi^\ast\R_{\tilde\MM}\to\R_T$ extends $\mu_1$ to $\mu$.  This
proves that every section of $(\tilde\MM^\sp,\R)$ comes locally from a
section of $\Modtilde_\alpha^\sp(A)$.  This finishes the proof.
\end{pf}

\subsection{The tangent complex}

Suppose $\LL=\LL^{\geq2}$ is a bundle of marked differential graded Lie algebras
on the smooth scheme (or algebraic space) $M$, and let $X\subset M$ be
its Maurer-Cartan locus. As a direct consequence of the second axiom
(Definition~\ref{bom}), the restriction of $(\LL,\delta)$ to $X$ is a
complex of sheaves of $\OO_X$-modules.  The derivative of the
marking $f:M\to \LL^2$ gives rise to an $\OO_X$-linear map
$T_M|_X\to\LL^2|_X$, and we obtain an augmented complex
$$\Theta^\bullet=\big[T_M|_X\longrightarrow\LL^2[1]|_X
\longrightarrow\LL^3[1]|_X\longrightarrow\ldots\big]$$  
by the first axiom. This complex $\Theta^\bullet$ of vector bundles on $X$ is
called the {\em tangent complex }of $(M,\LL)$. The shifts are applied
to $\LL|_X$ so that $T_M|_X$ will end up in degree $0$. 

By construction, the $0$-th cohomology sheaf is equal to the Zariski
tangent sheaf of $X$:
$$h^0(\Theta^\com)=T_X\,.$$

Next, we will recall how $\Theta$ governs deformation and obstruction theory.

\subsubsection{Deformation theory for small extensions}

Consider a pointed differential graded algebra $A\to\cc$,
concentrated in non-positive degrees.  Let $A'\to A$ be a small
extension of differential graded algebras: this means that the kernel
$I$ defined by the short exact sequence 
$$0\longrightarrow I\longrightarrow A'\longrightarrow
A\longrightarrow0\,,$$
and the kernel of the augmentation $A'\to\cc$ annihilate each
other. This implies that the $A'$-module structure on $I$ is induced
from the $\cc$-vector space structure on $I$ via the augmentation
$A'\to\cc$. For simplicity, assume that $I$ is concentrated in a specific
degree $-r\leq0$. (The classical case is the case $r=0$.) Denote by
$(T,\R_T)$ and $(T',\R_{T'})$ the affine 
differential graded schemes associated to $A$ and $A'$. 

We will consider a diagram
\begin{equation}\label{ext}
\vcenter{\xymatrix{\spec\cc\rto\drto\ar@/^1.5pc/[rrr]^P &
(T,\R_T)\rrto^{(\phi,\mu)}\dto & & (M,\R_M)\\ &
(T',\R_{T'})\ar@{.>}[urr]_{(\phi',\mu')}&&}}
\end{equation}
and ask for an obstruction to the existence of the dotted arrow.  If
a  dotted arrow exists, we will classify all possible dotted arrows up
to homotopy equivalence (see, for example, \cite{DerQuot} or
\cite{manetti}), for the definition of homotopy equivalence).

\begin{prop}\label{defob}
  There exists a naturally defined element $h\in
  H^1(P^\ast\Theta\otimes I)$, which vanishes if and only if a dotted
  arrow exists in~(\ref{ext}). If $h=0$, then the set of all dotted
  arrows in (\ref{ext}), up to homotopy equivalence, is in a natural
  way a principal homogeneous space for the abelian group
  $H^0(P^\ast\Theta\otimes I)$.
\end{prop}
\begin{pf}
The morphism $(\phi,\mu):(T,\R_T)\to(M,\R_M)$ is given by a morphism
of schemes $\phi:T\to M$ and a Maurer-Cartan element $\mu\in
MC(\phi^\ast\LL\otimes_{A^0}A)$. As $M$ is smooth, there is no
obstruction to the existence of $\phi'$, so let us choose $\phi'$. Now
consider the square zero extension of marked differential graded Lie
algebras
\begin{equation}\label{exts}
0\longrightarrow P^\ast\LL\otimes I\longrightarrow
{\phi'}^\ast\LL\otimes_{{A'}^0}A' \longrightarrow
\phi^\ast\LL\otimes_{A^0}A\longrightarrow0\,.
\end{equation}
We have a Maurer-Cartan element $\mu$ in the marked differential
graded Lie algebra on the right, which means that
$$f-\delta\mu+\sfrac{1}{2}[\mu,\mu]=0\,.$$
We lift $\mu$ at random to an element $\mu'$ of the marked differential graded Lie algebra
in the middle.  The obstruction $h$ is defined as
$$h=f'-\delta\mu'+\sfrac{1}{2}[\mu',\mu']\,,$$
which is an element of $P^\ast\LL\otimes I$, and moreover a 2-cocycle
in $P^\ast\LL\otimes I$, hence a 1-cocycle in $P^\ast\Theta\otimes
I$. The proof that vanishing of $h$ in cohomology is equivalent to the
existence of the dotted arrow distinguishes between the cases that
$r=0$ and $r<0$.  For $r<0$, we have $H^2(P^\ast\LL\otimes
I)=H^1(P^\ast\Theta\otimes I)$, and changing $\phi'$ while fixing
$\phi$ is impossible.  So the question is if there exists $z\in
P^\ast\LL\otimes I$ of degree $1$, such that $\mu'+z$ is a
Maurer-Cartan element in the middle of~(\ref{exts}). Such a $z$ will
exhibit $h$ as a coboundary (and conversely). For $r=0$, the element
$h\in H^1(P^\ast\Theta\otimes I)$ is the classical obstruction to the
existence of the dotted arrow in the diagram of classical schemes
$$\xymatrix{
\spec\cc\rto\drto & T\rrto\dto && X \ar@{^{(}->}[r]& M\\
& T'\ar@{.>}[urr]&}$$

Now assume that the obstruction vanishes. The difference between any
two Maurer-Cartan lifts of $\mu$ defines an element of
$H^0(P^\ast\Theta\otimes I)$.  One checks that this
difference is a coboundary if and only if the two lifts define
homotopy equivalent dotted arrows.
\end{pf}

\begin{cor}\label{cd}
For example, if $I=\cc[r]$, then the obstructions are contained in
$H^{r+1}(P^\ast\Theta)$ and the deformations are classified by
$H^r(P^\ast\Theta)$. 
\end{cor}

\subsubsection{Deformations of modules}

Let us examine the meaning of Proposition~\ref{defob} for the
differential graded algebraic space
$(\tilde\MM^\sp,\R)=\RModtilde_\alpha^\sp(A)$. 

So let the $\cc$-valued point $P:\spec \cc\to(\tilde\MM^\sp,\R)$ be represented
by the Maurer-Cartan element $\mu\in L$.

\begin{lem}
The complex $(L, d^{\mu})$ is precisely the graded normalized
Hochschild cochain complex with coefficients in $(\End V, \mu)$,
i.e. $\End V$ endowed with the structure of an $A$-$A$-bimodule from
$\mu$. 
\end{lem}
\begin{pf}
This is immediate. The normalized or reduced complex is defined, for
example, in 1.5.7 of~\cite{Loday}. 
\end{pf}

\begin{cor}
The complex $P^\ast\Theta$ is quasi-isomorphic to the augmented
graded Hochschild complex
$$\cc\longrightarrow \End_\gr V\longrightarrow\Hom_\gr(A,\End
V)\longrightarrow \Hom_\gr(A^{\otimes 2},\End V)\longrightarrow\ldots$$
\end{cor}
\begin{pf}
This follows immediately from the fact that the normalized Hochschild
complex is quasi-isomorphic to the Hochschild complex, by 1.5.7
of~\cite{Loday}. 
\end{pf}

\begin{cor}\label{exthom}
Suppose that
$P$ corresponds to the $A$-module $E=(V,\mu)$.
Then we have
$$H^i(P^\ast\Theta)=\begin{cases} \Ext_A^i(E,E)_\gr&\text{if
    $i>0$,}\\
\Hom_A(E,E)_\gr/\cc&\text{for $i=0$.}\end{cases}$$
The tangent complex $\Theta$ itself is quasi-isomorphic to the
augmented complex
$$\cc\longrightarrow R\sheafhom_A(\E,\E)_\gr\,,$$
where $\E$ is the universal family of graded $A$-modules  on
$\tilde\XX^\sp=\Modtilde_\alpha^\sp(A)$. 
\end{cor}
\begin{pf}
This is a consequence of the standard fact that (graded) Hochschild cohomology
computes (graded) extension spaces.  A proof in the ungraded case can
be found in~\cite{Weibel}, Lemma~1.9.1.
\end{pf}

\begin{cor}
In a situation given by a diagram such as~(\ref{ext}), assume that 
$I=\cc[r]$, as in Corollary~\ref{cd}.  
Then obstructions are contained
in $\Ext^{r+1}_A(E,E)_\gr$ and deformations are classified by
$\Ext^r_A(E,E)_\gr$ (or $\Hom_A(E,E)_\gr/\cc$, for $r=0$). 
\end{cor}

\section{Stability}

We will apply geometric invariant theory to the construction of the quotient of
$M=L^1$ by the gauge group $G$ as a quasi-projective scheme.

First note that since the scalars in $G$  act trivially,  no point of $L^1$
can be stable for the action of $G$. This prompts us to replace $G$ by
$\tilde G=G/\Delta$.  Second, notice that
the {\em canonical }one-parameter
subgroup $\lambda_0(t)=(t^p,\ldots,t^q)$ is central and acts by (see~(\ref{actij}))
$$(\lambda_0(t)\cdot \mu)_{ij}=t^{i-j}\mu_{ij}\,,$$
and hence destabilizes every element of $L^1$, as $i>j$ if
$\mu_{ij}\not=0$.  Thus the affine quotient $\spec \cc[L^{1}]^G$ is
trivial, equal to $\spec\cc$. 

In fact, the quotient of $L^1$ by $G$ classifies quiver
representations, for a certain quiver, and so the we are in the
situation worked out by King in~\cite{QuiverKing}. Our quiver has
$q-p+1$ vertices labelled $p,\ldots,q$, and for every pair of vertices
$i<j$, there are $\dim A_{j-i}$ arrows from $i$ to $j$. The vector
space $L^1=\Hom_\gr(\m,\End V)$ is denoted $\R(Q,\alpha)$
in~[ibid.], the group $G$ is denoted by $GL(\alpha)$. 

To linearize the action of $\tilde G$ on $L^1$, we choose a 
vector of integers
$$\theta=(\theta_p,\ldots,\theta_q)\,,\qquad\text{such that}\qquad
\sum_{i=p}^q\theta_i\alpha_i=0\,.$$ 
This defines the character $\chi_\theta:\tilde G\to\cc$ by
$$\chi_\theta(g)=\prod_{i=p}^q\det(g_i)^{\theta_i}\,,$$
which we use to linearize the action.

For a graded vector subspace $W\subset V$, define 
$$\theta(W)=\sum_{i=p}^q\theta_i\dim W_i\,.$$
Note that whether or not $\mu\in L^1$ satisfies the Maurer-Cartan
equation, it makes sense to speak of graded submodules $W\subset V$
with respect to $\mu$. 

\begin{prop}[King]
The point $\mu\in L^1$ is (semi)-stable for the action of $\tilde G$
linearized by $\chi_\theta$ if and only if for every proper graded $\mu$-submodule
$0<W<V$ we have $\theta(W)\mathrel{(\geq)}0$. (Here we use the usual convention
that to characterize stability, the strict inequality applies, and for
semi-stability the weak inequality is used.)
\end{prop}

Denote by $L^s$ and $L^{ss}$ the open subsets of $L^1$ of stable and
semi-stable points, respectively.  Similarly, denote by $X^s$ and
$X^{ss}$ the open subsets of stable and semi-stable points inside the
Maurer-Cartan subscheme $X\subset L^1$. 

The geometric invariant theory quotient of $L^1$ by $\tilde G$ is the
projective scheme
$$L^1\git\tilde G=\Proj \bigoplus_{n=0}^\infty\cc[L^1]^{G,\chi^n}\,,$$
where $\cc[L^1]^{G,\chi^n}=\{f:L^1\to\cc| f(gx)=\chi^n(g)f(x)\}$ is the
space of $\chi^n$-twisted invariants of $G$ in $\cc[L^1]$. The
quotient $L^1\git\tilde G$ is indeed projective, because
$\cc[L^1]^G=\cc$.

\begin{cor}
The scheme $L^s\git\tilde G$ is a quasi-projective smooth scheme
contained as an open subscheme in the algebraic space
$\tilde\MM^\sp$. It is a locally fine moduli space for equivalence
classes of stable quiver representations. In the coprime case, it is a
fine moduli space.

The scheme $L^1\git\tilde G=L^{ss}\git\tilde G$ is a projective scheme
containing $L^s\git\tilde G$ as an open subscheme.  Its points are in
one-to-one correspondence with $S$-equivalence classes of semi-stable
quiver representations. 
\end{cor}

\begin{cor}
The differential graded scheme $(L^s\git \tilde G,\R)$ is a
quasi-projective differential graded scheme, which represents the
sheaf associated to $\RModtilde_\alpha^s(A)$, the presheaf of
equivalence classes of families of stable graded $A$-modules. 

In the coprime case, $(L^s\git\tilde G,\R)$ represents
$\RModtilde_{\alpha}^s(A)$. 
\end{cor}

\begin{ex} \label{extreme-stability}
Maybe the most canonical of all characters is the one defined by
$\theta_p=-\dim V_q$, $\theta_q=\dim V_p$ and all other $\theta_i=0$.
We call it the {\em extremal }character. For this character, (semi)-stability reads
$$\dim W_p\dim V_q\mathrel{(\leq)}\dim W_q\dim V_p\,,$$
or, equivalently,
$$\frac{\dim W_p}{\dim W_q}\mathrel{(\leq)}\frac{\dim V_p}{\dim V_q}\,,$$
or
$$\frac{\dim W_p}{\dim V_p}\mathrel{(\leq)}\frac{\dim W_q}{\dim V_q}\,.$$
For example, stability implies that $V_p$ generates $V$ as an
$A$-module.
\end{ex}

\begin{defn}\label{module-stable}
We call the $[p,q]$-graded $A$-module $M$ (semi)-stable, if the
corresponding point $\mu$ in $L^1=\Hom_\gr\big(\m,\End(M)\big)$  is (semi)-stable with respect
to the linearization of $\tilde G$ given by the extremal character.
\end{defn}

\begin{ex}
Another canonical character is the determinant of the action of $G$ on
$L^1$. It has
$$\theta_i=\sum_{j<i}\dim A_{i-j}\dim V_j-\sum_{j>i}\dim
A_{j-i}\dim V_j\,,$$
and gives rise to the (semi)-stability condition
$$\sum_{i<j}\dim A_{j-i}\dim
W_i\dim V_j\mathrel{(\leq)}\sum_{i<j}\dim A_{j-i}\dim W_j\dim V_i\,.$$
\end{ex}

\section{Moduli of sheaves}

We will now assume that 
$A=\bigoplus_{n\geq0}\Gamma\big(Y,\OO(n)\big)$, for a connected projective
scheme $Y$.

\subsection{The adjoint of the truncation functor}

For a scheme $T$, we denote the projection $Y\times T\to T$ by
$\pi_T$. 

Let $T$ be a scheme and $\F_T$ a coherent sheaf on $Y\times T$. Then 
$$\Gamma_{[p,q]}\F_T=\bigoplus_{i=p}^q{\pi_T}_\ast\big(\F(i)\big)$$
is a graded sheaf of coherent $\OO_T$-modules with $A$-module
structure.

\newcommand{\ad}{{\mathscr{S}}}

\begin{prop}\label{adjoint}
The functor {\upshape
\begin{multline*}
\Gamma_{[p,q]}:(\text{coherent sheaves of $\OO_{Y\times T}$-modules})\\
\longrightarrow (\text{$[p,q]$-graded coherent sheaves of
  $A\otimes\OO_T$-modules})
\end{multline*}}has a left adjoint, which we shall denote by $\ad$. The functor $\ad$
commutes with arbitrary base change.
\end{prop}
\begin{pf}
First note that graded coherent $A\otimes\OO_T$-modules concentrated in the interval $[p,q]$
form an abelian category with kernels, cokernels, images and direct
sums constructed degree-wise, and that $\Gamma_{[p,q]}$ is an additive
functor, so that the statement makes sense. 

Then, by the claimed compatibility with base change, we may assume
that $T$ is affine, $T=\spec B$. 

Let $M$ be a graded $A\otimes B$-module, which
is concentrated in the interval $[p,q]$, and let
$$\xymatrix@1{
\bigoplus_j A(-m_j)\otimes B\rto & \bigoplus_i A(-n_i)\otimes B\rto & M\rto & 0}$$
be a presentation of $M$ (by graded homomorphisms) as a graded
$A\otimes B$-module.  Assume that all $n_i$ are in the interval
$[p,q]$.

Define $\ad M$ to be the cokernel in the diagram
of $\OO_{Y\times T}$-modules
\begin{equation}\label{defseq}
\xymatrix@1{
\bigoplus_{m_j\in[p,q]} \OO_{Y\times T}(-m_j)\rto & \bigoplus_{i}
\OO_{Y\times T}(-n_i)\rto & \ad M\rto & 0}\,,
\end{equation}
where the first sum extend only over those indices $j$, such that
$m_j$ is in 
the interval $[p,q]$. 
Let us prove that
$\ad M$ defined in this way satisfies
\begin{equation}\label{tws}
\Hom_{\OO_{Y\times T}}(\ad M,\fF)=\Hom_{A\otimes
  B}^\gr(M,\Gamma_{[p,q]}\fF)\,,
\end{equation}
for all $\OO_Y$-modules $\fF$. 
Given such $\fF$, consider the commutative diagram
$$\xymatrix{
\Hom_{\OO_{Y\times T}}\big(\bigoplus_i\OO_{Y\times T}(-n_i),\fF\big)\rto 
\ar@{=}[d]
& \Hom_{\OO_{Y\times T}}\big(\bigoplus_{m_j\in[p,q]}\OO_{Y\times
  T}(-m_j),\fF\big)
\ar@{=}[d]\\
\Hom_{A\otimes B}^\gr\big(\bigoplus_i A(-n_i)\otimes B,\Gamma_{[p,q]}\fF\big)\rto 
& \Hom_{A\otimes B}^\gr\big(\bigoplus_j A(-m_j)\otimes B,\Gamma_{[p,q]}\fF\big)\\
}$$
This diagram induces an equality of the kernels of the horizontal
maps, and these kernels are the two sides of (\ref{tws}), thus proving
Equation~(\ref{tws}). 

To prove that the adjoint functor $\ad$ commutes with base change,
consider a base change diagram
$$\xymatrix{
Y\times T'\rto\dto_v & T'\dto^u\\
Y\times T\rto & T}$$
and note that $\Gamma_{[p,q]}\circ v_\ast=u_\ast\circ\Gamma_{[p,q]}'$,
in obvious notation. It follows that for the adjoint functors we have
the equality $v^\ast\circ \ad=\ad'\circ u^\ast$, which is the claimed
compatibility with base change.
\end{pf}

\subsection{Open immersion}\label{regular}

Fix a numeric polynomial $\alpha(t)=a_0\binom{t}{0}+\ldots+a_k\binom{t}{k}$.

Let $\UU$ be a finite type open substack of the algebraic stack of
coherent sheaves on $Y$ with Hilbert polynomial $\alpha(t)$. 

For $q>p>0$ let $\Mod_\alpha^{[p,q]}(A)$ be the algebraic stack of $[p,q]$-graded
$A$-modules of dimension $\alpha|_{[p,q]}$. Recall that
$\Mod_\alpha^{[p,q]}(A)=[MC(L)/G]$, in the notation of Section~\ref{RAct}.

\begin{prop}
Given $\UU$, there exists $p$, such that for all $q>p$ the functor $\Gamma_{[p,q]}$
defines a morphism of algebraic stacks
$$\Gamma_{[p,q]}:\UU\longrightarrow\Mod_\alpha^{[p,q]}(A)\,.$$
If $q$ is sufficiently large, then $\Gamma_{[p,q]}$ is a monomorphism
of stacks.
\end{prop}
\begin{pf}
Let $p$ be large enough such that every sheaf in $\UU$ is Castelnuovo-Mumford
$p$-regular. Then, for every $i\geq p$, the sheaf $\pi_\ast \F(i)$ is
locally free of rank $\alpha(i)$ on $T$. Hence $\Gamma_{[p,q]}\F$ is a
$\Mod^{[p,q]}_\alpha(A)$-family over $T$, and we have the required
morphism of stacks. 

Now let, in addition, $q$ be large enough for $\OO_Y(q-p)$ to be Castelnuovo-Mumford
regular. Then $\Gamma_{[p,q]}$ is a monomorphism of stacks, because
for every family of $p$-regular sheaves $\F$, the adjunction map
$\ad(\Gamma_{[p,q]}\F)\to \F$ is an isomorphism. See \cite{Alv-King},
Theorem~3.4 and Proposition~4.1,
for a proof of a similar statement. In our context, we may proceed as
follows:

First note that we may assume that the parameter scheme $T$ is affine,
$T=\spec B$, as in the proof of Proposition~\ref{adjoint}.

Let $V=\Gamma(Y,\fF(p))$, and $\gG$ the kernel in:
$$\xymatrix@1{
0\rto & \gG\rto & V\otimes_B \OO_Y(-p)\rto & \fF\rto & 0}\,.$$
Then the fact that $\OO_Y(q-p)$ is regular implies that $\gG$ is
$q$-regular. (Lemma~3.3 in~\cite{Alv-King}.) Let
$W=\Gamma\big(Y,\gG(q)\big)$, so that we have a 
surjection $W\otimes_B \OO_Y(-q)\twoheadrightarrow\gG$, and a
presentation of $\fF$:
$$\xymatrix@1{
W\otimes_B \OO_Y(-q)\rto & V\otimes_B \OO_Y(-p)\rto & \fF \rto& 0}$$
We remark that $q$-regularity of $\gG$ implies that this sequence
stays exact after twisting by $\OO_Y(i)$ 
and taking global sections, for all $i\geq q$. Thus the
sequence of graded $A\otimes B$-modules 
\begin{equation}\label{pmr}\nonumber
\xymatrix@1{
W\otimes A(-q)\rto & V\otimes A(-p)\rto & \Gamma_{\geq p}\fF\rto &
0}
\end{equation}
is exact in degrees $\geq q$. 
We can construct from this a presentation of $\Gamma_{\geq p}\fF$, by
adding some relations whose degrees are between $p$ and $q$.  Then
we can turn this presentation of $\Gamma_{\geq p}\fF$ into a presentation
of $\Gamma_{[p,q]}\fF$ by adding relations in degrees larger than $q$.  These extra
relations in degrees larger than $q$ are ignored when
constructing $\ad(\Gamma_{[p,q]}\fF)$, see the proof of
Proposition~\ref{adjoint}. The extra relations of degree between $p$ and
$q$ do not affect the cokernel in~(\ref{defseq}). We conclude that we
have a presentation 
$$\xymatrix@1{
W\otimes_B \OO_Y(-q)\rto & V\otimes_B \OO_Y(-p)\rto &
\ad(\Gamma_{[p,q]}\fF)  \rto& 0}\,.$$ 
This proves that $\ad(\Gamma_{[p,q]}\fF)=\fF$.
\end{pf}

\begin{prop}
For $q\gg p\gg0$ the morphism
$\Gamma_{[p,q]}:\UU\to\Mod_\alpha^{[p,q]}(A)$ is \'etale.
\end{prop}
\begin{pf}
Let $A'\to A\to\cc$ be a small extension of pointed $\cc$-algebras (not
differential graded). Let $T=\spec A$ and $T'=\spec A'$. Consider a
2-commutative diagram
$$\xymatrix{
T\rrto^{\F}\ar@{^(->}[d] && \UU\dto^{\Gamma_{[p,q]}}\\
T'\rrto^{M'}\ar@{.>}[urr] && \Mod_\alpha^{[p,q]}(A)}$$
of solid arrows.  We have to prove that the dotted arrow exists,
uniquely, up to a unique 2-isomorphism. This follows from standard
deformation-obstruction theory.  We need that $\Gamma_{[p,q]}$ induces
a bijection on deformation spaces  and an injection on obstruction
spaces (associated to the above diagram). It is
well-known, that deformations of 
$\F$ are classified by $\Ext^1_{\OO_Y}(\F,\F)$ and obstructions are
contained in $\Ext^2_{\OO_Y}(\F,\F)$.  We saw in
Corollary~\ref{exthom}, that deformations of $M'|_T$ are classified by
$\Ext^1_A(M,M)_\gr$ and obstructions are contained in
$\Ext^2_A(M,M)_\gr$, where $M=\Gamma_{[p,q]}(\F)$. It is proved in
\cite{DerQuot}, (4.3.3.a) and (4.3.4), that for fixed $i$, there exist
$q\gg p\gg0$, such that
$\Ext^i_{\OO_Y}(\F,\F)=\Ext^i_A(\Gamma_{[p,q]}\F,\Gamma_{[p,q]}\F)_\gr$. (Note
that the assumption in \cite{DerQuot}, that $Y$ be smooth is not used
for this result.  It is only used to exchange quantifiers: namely to
get uniform $p$ and $q$, which work for {\em all }$i\geq0$.)
\end{pf}

\begin{cor}
For $q\gg p\gg0$ the morphism
$\Gamma_{[p,q]}:\UU\to\Mod_\alpha^{[p,q]}(A)$ is an open immersion.
\end{cor}

\subsection{Stable sheaves}

\newcommand{\ssgeq}{{> \atop {\scriptscriptstyle (}-{\scriptscriptstyle )} }}
\newcommand{\ssleq}{{< \atop {\scriptscriptstyle (}-{\scriptscriptstyle )} }}

Let $Y$ be a connected projective
scheme.
We denote the Hilbert polynomial of a coherent sheaf $\fF$ on $Y$ by
$h(\fF,t)=h(\F)$.

For our purposes, the following characterization of stability of $\F$
is most useful.

\begin{defn}
The sheaf $\F$ is {\em (semi)-stable }if and only if for every proper
subsheaf $0<\F'<\F$ we have
$$\frac{h(\F',p)}{h(\F,p)}\mathrel{(\leq)}
\frac{h(\F',q)}{h(\F,q)}\qquad\text{for $q\gg p\gg 0$}\,.$$
(As usual, this means the strict inequality for `stable' and the weak
inequality for `semi-stable'.)

The condition needs only to be checked for {\em saturated
}subsheaves. (A subsheaf is saturated if the corresponding quotient is
pure of the same dimension as $\F$.)
\end{defn}

\begin{rmk}
We can say, informally, that the limiting slope of the quotient of
Hilbert polynomials $\frac{h(\F')}{h(\F)}$ is $\mathrel{(\geq)}0$, for
all proper saturated subsheaves. 
\end{rmk}

This stability condition looks very similar to the condition given by
the extremal character for $A$-modules, see
Example~\ref{extreme-stability}, but to relate the two notions is not 
completely trivial.

\begin{thm}\label{prescribe}
Given $\UU$, it is possible to choose $q\gg p\gg0$ in such a way that the following
holds:
if $\F$ is a $\UU$-sheaf, then $\fF$ is a
(semi)-stable sheaf if and only if 
$M=\Gamma_{[p,q]}\fF$ is a (semi)-stable graded $A$-module
(Definition~\ref{module-stable}). 
\end{thm}
\begin{pf}
By Grothendieck's lemma (see \cite{HuybrechtsLehn}, Lemma~1.7.9), the
family $\UU'$ of all saturated destabilizing subsheaves of all sheaves
in $\UU$ is bounded. We choose $p$ large enough to ensure that all
sheaves in $\UU$ and $\UU'$ are
$p$-regular. Note that the sheaves in $\UU'$ have only finitely
many Hilbert polynomials.  So we can choose $q\gg p\gg 0$ in such a
way that the limiting slope of all quotients of all Hilbert
polynomials involved  is measured correctly by $p$ and $q$.

Additionally, we choose $p$ and $q$ sufficiently large as explicated
in~\cite{Alv-King}. (This choice is only needed for the `converse', below.)

Let us first suppose that $M$ is (semi)-stable, and prove that $\F$ is
(semi)-stable. So let $0\subsetneq\F'\subsetneq \F$ be a saturated
subsheaf. We wish to prove, of course, that $\F'$ does not violate
(semi)-stability of $\F$. So let us assume it does. Then by our
choices, both $\F'$ and $\F$ are $p$-regular.

Since
$\Gamma_{[p,q]}$ is left exact, we get a graded submodule
$$M'=\Gamma_{[p,q]}\F'\hooklongrightarrow\Gamma_{[p,q]}\F\,.$$
Moreover, $0\subsetneq M'\subsetneq M$, as
$\F'=\ad M'$, because $\F'$ is $p$-regular.
Since $M=\Gamma_{[p,q]}\F$ is (semi)-stable, we know that
$$\frac{\dim\Gamma(Y,\F'(p))}{\dim\Gamma(Y,\F(p))}\mathrel{(\leq)}
\frac{\dim\Gamma(Y,\F'(q))}{\dim\Gamma(Y,\F(q))}\,.$$
By $p$-regularity, this implies that
$$\frac{h(\F',p)}{h(\F,p)}\mathrel{(\leq)}
\frac{h(\F',q)}{h(\F,q)}\,,$$
and so $\F'$ does {\em not }violate (semi)-stability, a
contradiction.

Conversely, assume that $\F$ is (semi)-stable. If $0<M'<M$ is a
(semi)-stability violating submodule, then $(M'_p,M'_q)\subset
(M_p,M_q)$ is a Kronecker submodule in the sense of
\cite{Alv-King}. To prove that $(M'_p,M'_q)\not=(0,0)$, note that
$\Gamma_{[p,q]}\F$ does not have any non-trivial submodules which
  vanish in the top degree $q$. (This is an elementary fact
  about sheaves on projective schemes.) To prove that
  $(M'_p,M'_q)\not=(M_p,M_q)$, note that $\Gamma_{[p,q]}\F$ is
  generated in the lowest degree $p$, by $p$-regularity of $\F$. 

Thus, applying Theorem~5.10 of [ibid.], we see that $M'$ does {\em not
}violate (semi)-stability, a contradiction. 
\end{pf}

\subsection{Moduli of Sheaves}

Let $\alpha(t)$ be a Hilbert polynomial. Let $\UU^{ss}$ be the bounded
family of all semi-stable sheaves with Hilbert polynomial $\alpha(t)$.
Choose $q\gg p\gg0$ as prescribed by Theorem~\ref{prescribe} for
$\UU^{ss}$. 

Let $U^s\subset U^{ss}$ be the moduli spaces of stable
(resp. semistable) sheaves on $Y$ with Hilbert polynomial
$\alpha(t)$. 

\begin{cor}
We have a commutative diagram of open immersions of schemes.
$$\xymatrix{
U^{ss}\rrto^-{\Gamma_{[p,q]}} && L^{ss}\git \tilde G\\
U^s\uto\rrto^-{\Gamma_{[p,q]}} &&
\Modtilde^s_{\alpha|_{[p,q]}}(A)\uto}$$
The two schemes in the top row are projective. Hence, $U^{ss}$ is a
union of connected components of
$L^{ss}\git\tilde G$.

In the case where
$\alpha$ is indivisible, we have $U^s=U^{ss}$, and so $U^s$ is a union
of components of $\Modtilde^s_{\alpha|_{[p,q]}}$, via the functor
$\Gamma_{[p,q]}$.
\end{cor}

\begin{rmk}
Assume we are in the indivisible case. Then $U^s\subset L^s\git\tilde
G$ is a closed subscheme of the smooth scheme $L^s\git\tilde G$, cut
out by the descended Maurer-Cartan equation $dx+\frac{1}{2}[x,x]=0$.
This gives rather explicit equations for $U^s$ inside a smooth
scheme. Note that we do not prove that $L^s\git\tilde G$ is
projective, in the indivisible case. 
\end{rmk}

\subsection{An amplification}

By using three integers $q\gg p'\gg p\gg0$, we can describe the image
of $\Gamma_{[p,q]}:\UU^{s}\to\Mod_{\alpha|_{[p,q]}}(A)$ explicitly. 

We denote by $\Mod_{\alpha|_{[p,q]}}(A)'\subset\Mod_{\alpha|_{[p,q]}}(A)$
the open substack of graded $A$-modules which are generated in degree
$p$.

\begin{thm}
Let $\UU$ be, as above, a bounded open family of sheaves on $Y$. Then
for $q\gg p'\gg p\gg 0$, the functor $\Gamma_{[p,q]}$ induces an open
immersion
$$\Gamma_{[p,q]}:\UU\longrightarrow \Mod_{\alpha|_{[p,q]}}(A)'$$
and the image of $\UU^{s}$ ($\UU^{ss}$) is equal to the locus of
modules 
whose truncation into the interval $[p',q]$ is (semi)-stable.
\end{thm}
\begin{pf}
The first claim is clear: $p$-regularity of $\F$ implies that
$\Gamma_{[p,q]}\F$ is generated in degree $p$.

The fact that $\UU^{s}$ ($\UU^{ss}$) is contained in the
$[p',q]$-(semi)-stable locus follows from Theorem~\ref{prescribe}. 

Let $M$ be an $A$-module concentrated in degrees $[p,q]$, generated in
degree $p$, and of dimension $\alpha|_{[p,q]}$. Then we will use
Gotzmann persistence to prove that $\F=\ad(M)$ has
Hilbert polynomial $\alpha$, and we will make sure that all 
$\ad(M)$ obtained in this way are $p'$-regular. This will imply that 
$M_{[p',q]}=\Gamma_{[p',q]}\F$, and we can again apply
Theorem~\ref{prescribe} to deduce that if $M_{[p',q]}$ is
(semi)-stable, then $\F$ is (semi)-stable. 

We briefly recall the persistence theorem (see \cite{gotz} and
\cite{gash}, especially Theorem~4.2 in \cite{gash}). First, for integers $a\geq0$ and $t\geq1$, there exist
unique integers $m_t>m_{t-1}>\ldots>m_1\geq0$, such that
$a=\sum_{i=1}^t \binom{m_i}{i}$. Then one defines $a^{\langle
  t\rangle}=\sum_{i=1}^t\binom{m_i+1}{i+1}$. One significance of this
definition is the following: if $\E$ is a coherent sheaf of
$\OO_Y$-modules, such that $\E(p)$ is globally generated, and if $h(t)$
is the Hilbert polynomial of $\E$, then $h(t+1)=h(t)^{\langle
  t-p\rangle}$, for $t\gg0$. The persistence theorem says the
following:

Suppose $A$ is a graded $\cc$-algebra, generated in degree 1, with
relations in degree $\leq r$, for an integer $r\geq 1$. Let $M$ be a
graded $A$-module and $G$ a finite dimensional graded $\cc$-vector
space, such that the following sequence of graded $A$-modules is
exact:
$$0\longrightarrow K\longrightarrow A\otimes_\cc G\longrightarrow
M\longrightarrow 0\,.$$
\begin{items}
\item (Macaulay bound) If $\deg G\leq p$, then $\dim M_{d+1}\leq(\dim M_d)^{\langle
    d-p\rangle}$, for all $d\geq p+1$.  Moreover, there exists a $d$, such that $\dim M_{d'+1}=(\dim
  M_{d'})^{\langle d'-p\rangle}$, for all $d'\geq d$. .   
\item (Persistence) If, in addition,  $K$ is generated in degree $\leq r'$, where $r'\geq p+r$, and if
  $\dim M_{d+1}=(\dim M_d)^{\langle d-p\rangle}$, for some $d\geq r'$, then $\dim
  M_{d'+1}=(\dim M_{d'})^{\langle d'-p\rangle}$, for all $d'\geq d$.
\end{items}
We may assume, that $\alpha(t+1)=\alpha(t)^{\langle t-p\rangle}$, for
all $t\geq p$.

Now let $M$ be an $A$-module in $[p,q]$ of  dimension
$\alpha|_{[p,q]}$, which is generated in degree $p$. We have the exact sequence
$$0\longrightarrow K\longrightarrow A_{[0,q-p]}\otimes
M_p\longrightarrow M\longrightarrow 0\,,$$
where the kernel $K$ exists (at most) in degrees $[p+1,q]$. 
Let $\tilde K\subset A\otimes M_p$ be the submodule generated by $K$,
and let $\tilde M$ be the quotient
$$0\longrightarrow \tilde K\longrightarrow A\otimes
M_p\longrightarrow \tilde M\longrightarrow 0\,.$$
Thus $\tilde K$ is generated in degree $\leq q$. 

Or first claim is that $\tilde K$ is actually generated in degree
$p+1$. We will do
this by descending induction. So suppose $\tilde K$ is generated in
degree $\leq r'$, for $p+1<r'\leq q$, but not in degree $\leq r'-1$. Then
let $\tilde K'< \tilde K$ be the submodule generated by the degree
$\leq r'-1$ part of $\tilde K$. Let $\tilde M'=(A\otimes M_p)/\tilde
K'$ be the quotient. Then we have
$$(\dim \tilde M'_{r'-1})^{\langle r'-1-p \rangle}\geq \tilde
M'_{r'}>\tilde M_{r'}=(\tilde M_{r'-1})^{\langle r'-1-p\rangle}\,,$$
which implies $\dim \tilde M'_{r'-1}>\dim \tilde M_{r'-1}$, which is
absurd, as these two spaces are equal.  Thus $\tilde K$ is, indeed,
generated in degrees $\leq r'-1$, and we conclude that it is, in fact,
generated in degree $p+1$. 

Now, the persistence
theorem implies that $\dim \tilde M_{t+1}=\dim\tilde
M_{t}^{\langle t-p\rangle}$, for all $t> p+r$. As $\ad(M)$ is the
sheaf associated to $\tilde M$, this implies that the Hilbert
polynomial of $\ad(M)$ is equal to $\alpha$, as claimed. 

We remark that the family of all $A$-modules generated in degree $p$
by $\alpha(p)$ elements, whose relations are in degree $p+1$, is
bounded. Therefore, we can choose $p'>p$ in such a way that all sheaves
associated to such modules are $p'$-regular.  This will imply that  all $\ad(M)$ obtained from 
$\Mod_{\alpha|_{[p,q]}}(A)'$ are  $p'$-regular.

It remains to prove that a suitable choice of $p'$ will assure that
the truncation of $M$ into the interval $[p',q]$ is equal to
$\Gamma_{[p',q]}\F$, where $\F=\ad(M)$. 

Now, the canonical map $\tilde M\to \Gamma_{\geq p}\F$ is an
isomorphism in sufficiently high degree. But as the
family of all $\tilde M$ which occur is bounded, there exists a
uniform $p'$ which will assure that $\tilde M|_{\geq p'}\to\Gamma_{\geq p'}\F$ is an
isomorphism. This finishes the proof of the last remaining fact that
$M|_{[p',q]}=\Gamma_{[p',q]}\ad(M)$. 
\end{pf}

\begin{cor}
We have
$$U^{ss}=\Modtilde^{{\scriptscriptstyle{[p',q]}}-ss}_{\alpha|_{[p,q]}}(A)'\,,$$
and 
$$U^{s}=\Modtilde^{\scriptscriptstyle{[p',q]}-s}_{\alpha|_{[p,q]}}(A)'\,,$$
in obvious notation.
In the indivisible case, all four schemes are equal.
\end{cor}

\begin{rmk}
If an $A$-module in $[p,q]$ is stable (not just semi-stable),
then it is generated in degree $p$.
Thus $U^s$ can also be described as the scheme of modules in the
interval $[p,q]$, of dimension $\alpha|_{[p,q]}$, which are stable,
and whose truncation into the interval $[p',q]$ is also stable. 
\end{rmk}

\section{Derived moduli of sheaves}\label{4.0}

Finally, we will construct the differential graded moduli scheme of
stable sheaves on the projective variety $Y$.  From now on, we have to
assume that $Y$ is smooth. Let $\alpha(t)$ be a numerical polynomial,
and $p\gg0$. For simplicity, let us assume that $\alpha(t)$ is
primitive. 

\begin{defn}
A {\em family of coherent sheaves }on $Y$, of Hilbert polynomial $\alpha(t)$, parametrized by the differential
graded scheme $(T,\R_T)$, is a pair $(E,\mu)$, where $E$ is a graded quasi-coherent sheaf
$$E=\bigoplus_{i\geq p} E^i$$
on $T$, and each $E_i$ is a vector bundle, of rank $\rk
E_i=\alpha(i)$. Moreover, $\mu$ is a `unitary' 
Maurer-Cartan element in the differential graded Lie algebra
$$\Gamma\big(T,\sheafhom_\gr(A^{\otimes\geq
  1},\sheafend_{\OO_T}E)\otimes_{\OO_T}\R_T\big),$$
in other words a graded unitary $\R_T$-linear $A_\infty$-action of $A\otimes
\R_T$ on $E\otimes_{\OO_T}\R_T$.

We denote the functor of equivalence classes of simple  such families by
$\RModtilde_\alpha^{\sp}(\OO_Y)$. 

If for every point $P:\spec\cc\to T$, the associated coherent sheaf on
$Y$ is (semi)-stable, then the family $(E,\mu)$ is a {\em
  (semi)-stable }family.

\end{defn}

\begin{lem}
We have
$$\RModtilde_\alpha^\sp(\OO_Y)=\projectlim_{q\gg
  p}\RModtilde_{\alpha|_{[p,q]}}^\sp(A)\,,$$
as set-valued presheaves on the category of differential graded
schemes. 
\end{lem}
\begin{pf}
{\em Obvious}.
\end{pf}

\begin{cor}
The functor $\RModtilde_\alpha^\sp(\OO_Y)$ is represented by the projective
limit of differential graded algebraic spaces
$$\RModtilde_\alpha^\sp(A)=\projectlim_{q\gg
  p}\RModtilde_{\alpha|_{[p,q]}}^\sp(A)\,.$$
\end{cor}

\begin{prop}
The projective limit
$$\projectlim_{q\gg
  p}\RModtilde_{\alpha|_{[p,q]}}^\sp(A)$$
stabilizes, as far as quasi-isomorphism is concerned.
\end{prop}
\begin{pf}
Here we use that $Y$ is smooth, to deduce that 
$$\Ext_{\OO_Y}^i(E,E)=\Ext_{A}^i(\Gamma_{[p,q]}E,\Gamma_{[p,q]}E)\,,$$
for $q\gg p$. Then we use the fact that if $\pi_0$ agrees, and tangent
complex cohomologies agree, then a morphism of differential graded
schemes is a quasi-isomorphism.
\end{pf}

\begin{cor}
If $q\gg p$, then $\RModtilde_\alpha^\sp(\OO_Y)$ is quasi-isomorphic to
$\RModtilde_{\alpha|_{[p,q]}}^\sp(A)$. Moreover,
$\RModtilde_\alpha^s(\OO_Y)$ is an open and closed differential graded
subscheme of $\RModtilde_{\alpha|_{[p,q]}}^\sp(A)$. 
\end{cor}

\bibliographystyle{plain}

\begin{thebibliography}{10}

\bibitem{Alv-King}
L.~{\'A}lvarez-C{\'o}nsul and A.~King.
\newblock A functorial construction of moduli of sheaves.
\newblock {\em Invent. Math.}, 168(3):613--666, 2007.

\bibitem{BF}
K.~Behrend and B.~Fantechi.
\newblock The intrinsic normal cone.
\newblock {\em Invent. Math.}, 128(1):45--88, 1997.

\bibitem{DerQuot}
I.~Ciocan-{F}ontanine and M.~Kapranov.
\newblock Derived {Q}uot schemes.
\newblock {\em Ann. Sci. \'Ecole Norm. Sup. (4)}, 34(3):403--440, 2001.

\bibitem{vfcdgm}
I.~Ciocan-Fontanine and M.~Kapranov.
\newblock Virtual fundamental classes via dg-manifolds.
\newblock {\em Geom. Topol.}, 13(3):1779--1804, 2009.

\bibitem{gash}
V.~Gasharov.
\newblock Extremal properties of {H}ilbert functions.
\newblock {\em Illinois J. Math.}, 41(4):612--629, 1997.

\bibitem{gotz}
G.~Gotzmann.
\newblock Eine {B}edingung f\"ur die {F}lachheit und das {H}ilbertpolynom eines
  graduierten {R}inges.
\newblock {\em Math. Z.}, 158(1):61--70, 1978.

\bibitem{HuybrechtsLehn}
D.~Huybrechts and M.~Lehn.
\newblock {\em The geometry of moduli spaces of sheaves}.
\newblock Aspects of Mathematics, E31. Friedr. Vieweg \& Sohn, Braunschweig,
  1997.

\bibitem{QuiverKing}
A.~D. King.
\newblock Moduli of representations of finite-dimensional algebras.
\newblock {\em Quart. J. Math. Oxford Ser. (2)}, 45(180):515--530, 1994.

\bibitem{Loday}
J.-L. Loday.
\newblock {\em Cyclic homology}, volume 301 of {\em Grundlehren der
  Mathematischen Wissenschaften}.
\newblock Springer-Verlag, Berlin, 1992.

\bibitem{lurie}
J.~Lurie.
\newblock {\em Higher topos theory}, volume 170 of {\em Annals of Mathematics
  Studies}.
\newblock Princeton University Press, Princeton, NJ, 2009.

\bibitem{manetti}
M.~Manetti.
\newblock Deformation theory via differential graded {L}ie algebras.
\newblock In {\em Algebraic {G}eometry {S}eminars, 1998--1999 ({P}isa)}, pages
  21--48. Scuola Norm. Sup., Pisa, 1999.
\newblock See also arXiv:math.AG/0507284.

\bibitem{simpson}
C.~T. Simpson.
\newblock Moduli of representations of the fundamental group of a smooth
  projective variety. {I}.
\newblock {\em Inst. Hautes \'Etudes Sci. Publ. Math.}, 79:47--129, 1994.

\bibitem{HolCas}
R.~P. Thomas.
\newblock A holomorphic {C}asson invariant for {C}alabi-{Y}au 3-folds, and
  bundles on {$K3$} fibrations.
\newblock {\em J. Differential Geom.}, 54(2):367--438, 2000.

\bibitem{ModuliDGCats}
B.~To{\"e}n and M.~Vaqui{\'e}.
\newblock Moduli of objects in dg-categories.
\newblock {\em Ann. Sci. \'Ecole Norm. Sup. (4)}, 40(3):387--444, 2007.

\bibitem{BraveNewPaper}
B.~To{\"e}n and G.~Vezzosi.
\newblock From {HAG} to {DAG}: derived moduli stacks.
\newblock In {\em Axiomatic, enriched and motivic homotopy theory}, volume 131
  of {\em NATO Sci. Ser. II Math. Phys. Chem.}, pages 173--216. Kluwer Acad.
  Publ., Dordrecht, 2004.

\bibitem{HOGI}
B.~To{\"e}n and G.~Vezzosi.
\newblock Homotopical algebraic geometry. {I}. {T}opos theory.
\newblock {\em Adv. Math.}, 193(2):257--372, 2005.

\bibitem{HOGII}
B.~To{\"e}n and G.~Vezzosi.
\newblock Homotopical algebraic geometry. {II}. {G}eometric stacks and
  applications.
\newblock {\em Mem. Amer. Math. Soc.}, 193(902):x+224, 2008.

\bibitem{Weibel}
C.~A. Weibel.
\newblock {\em An introduction to homological algebra}, volume~38 of {\em
  Cambridge Studies in Advanced Mathematics}.
\newblock Cambridge University Press, Cambridge, 1994.

\end{thebibliography}

\end{document}